\documentclass[11pt, twoside]{article}

\usepackage[english]{babel}
\usepackage[T1]{fontenc} 
\usepackage[latin1]{inputenc} 

\usepackage{amsmath,amsthm,amscd,latexsym, mathtools, eucal,epsfig,enumitem,bbold, color} 
\usepackage{pstricks,pstricks-add,pst-math,pst-xkey,amssymb,stmaryrd,ulem}  
\usepackage{amsmath,amsthm,amscd,latexsym, mathtools, eucal,epsfig,enumitem,bbold, color} 
\usepackage[bookmarks={true},bookmarksopen={true}]{hyperref} 
\let\norm\undefined 
\DeclarePairedDelimiter\norm{\|}{\|} 

\newtheorem{theorem}{Theorem}
\newtheorem{proposition}[theorem]{Proposition}
\newtheorem{lemma}[theorem]{Lemma}

\theoremstyle{definition}
\newtheorem{definition}[theorem]{Definition}
\newtheorem{notation}[theorem]{Notation} 
\newtheorem{remark}[theorem]{Remark}

\theoremstyle{theorem}
\newtheorem*{mainResultI}{Main Result}
\newtheorem*{corollary*}{Corollary}

\newtheorem*{remark*}{Remark}
\theoremstyle{definition}
\newtheorem*{definition*}{Definition}
\begin{document}

\title{Ordinary differential equations defined by a trigonometric polynomial field:   Behavior  of   the solutions  }

\author{W. Oukil\\
\small\text{Laboratory of Dynamical Systems, Faculty of Mathematics. }\\
\small\text{University of Science and Technology Houari Boumediene,}\\
\small\text{ BP 32 EL ALIA 16111 Bab Ezzouar, Algiers, Algeria}. }

\date{\today}

\maketitle

\begin{abstract}
We consider the ordinary differential equations defined by a trigonometric polynomial field, we prove   that any solution $x$  admits a {\it rotation vector} $\rho\in \mathbb{R}^n$. More precisely, the function $t\mapsto x(t)-\rho t$ is bounded on time and it is a {\it weak almost periodic} function of {\it slope} $\rho$.
\end{abstract}
\\




 \section{Introduction}
In this article,  we study the asymptotic behavior  of   solutions for    ordinary differential equations (ODE) defined by a trigonometric polynomial field. The idea comes from the scalar case, where in this case  H. Poincar\'e defined the {\it rotation number} for circle homeomorphisms  \cite{Poincare}. The simple example is a   scalar differential equation 
\[
\dot{x}=f(x),\quad x(0)\in\mathbb{R},\quad t\in\mathbb{R},
\]
where $f:\mathbb{R}\to\mathbb{R}$ is  lipschitz,  1-periodic and $t\mapsto x(t)$ is the  state of the system. There exists  a rotation number  $\lambda\in\mathbb{R}$ for which the function $t\mapsto x(t)-\lambda t$ is   bounded (periodic). We know that any non-autonomous system can be written as  an autonomous system. Our result is a generalization of this asymptotic behavior to any dimension. In this case, $\lambda$ is a vector and called a {\it rotation vector} or  {\it rotation set} as it is defined in \cite{Michal}.  Under some assumptions  of stability [\cite{Saito1971_01}, \cite{Fink1974_01}]  proved the existence of the rotation vector.   Some biological works   use the ODE defined by a trigonometric polynomial field and study the rotation vector components  as in [\cite{AriaratnamStrogatz}, \cite{Ha}, \cite{kuramoto1}, \cite{WinfreeModel}]. Our contribution to this  biological works has two key points, the mathematical proof of   existence of the rotation vector and   the study of the behavior of solutions. 
\section{Definition and Main result}  
\noindent We study in this article the following system
\begin{equation*}
\dot{x}=  f(x ),\quad t\in\mathbb{R},\quad  x(0) =  x_0\in \mathbb{R}^n,
\end{equation*}
where $t\mapsto x(t):=(x_j(t))_{j=1}^n$ is the  state of the system and $f:  \mathbb{R}^n\to \mathbb{R}^n$  is  a trigonometric polynomial  in the following sense
\begin{definition}\label{def:tpf}[Trigonometric polynomial function]
A function $g: \mathbb{R}^n \to \mathbb{R}$ is called {\it a trigonometric polynomial} if there exists a finite sequence $(c_p)_{p\in 2\pi\mathbb{Z}^n}\subset \mathbb{C}$  such that
\[
\forall x\in\mathbb{R}^n:\quad g(x) =\sum_{p\in 2\pi\mathbb{Z}^n }c_p\exp( i\langle x,p\rangle ),
\]
where $\langle .,.\rangle$ is the usual scalar product on $\mathbb{R}^n$. A function $g: \mathbb{R}^n \to \mathbb{R}^n$ is {\it a trigonometric polynomial} if each  component is a trigonometric polynomial function.
\end{definition}
To formulate the  Mains results   let us introduce the following definitions.  We use the usual  norm $\norm{y}:=\max_{1\le j\le n}\norm{y_j}$ for every $y:=(y_j)_{j=1}^n\in\mathbb{C}^n$ .  
\begin{definition}\label{def:almostperiofuncmodulo}[Rotation vector]
Let $\lambda\in\mathbb{R}^n$ and $\phi:\mathbb{R}\to\mathbb{R}^{n}$ be a function. We say that $\phi$ admits $\lambda$ as the {\it rotation vector } if 
\[
\sup_{t\in\mathbb{R}}\norm{\phi(t)-\lambda t}<\infty. 
\]
\end{definition}
For more information about the behavior of solutions, we introduce the following definitions.
\begin{definition}\label{def:periomodulo}[Periodic   modulo $\mathbb{Z}^n$  function]
A function   $g: \mathbb{R}^n \to \mathbb{R}^n$   is called {\it  periodic modulo $\mathbb{Z}^n$}, if 
\[
g(z_1+k_1,\ldots, z_n+k_n)=g(z_1 \ldots, z_n),\quad \forall (k_j)_{j=1}^n \in \mathbb{Z}^n,\ \forall (z_j)_{j=1}^n\in \mathbb{R}^n.
\]
\end{definition}
 \begin{definition}\label{def:almostperiofunc}[Weakly almost-periodic function]
Let be $r\in \mathbb{R}^n$. A function $h: \mathbb{R}  \to \mathbb{R}^n$ is   {\it weakly almost periodic of slope $r$}  if it is $C^\infty$ and if there exists a uniformly bounded sequence  for the sup-norm  of  $C^\infty$  functions $(g_k: \mathbb{R}^n \to \mathbb{R}^n)_{k\in\mathbb{N}}$ that are periodic modulo $\mathbb{Z}^n$  and there exists  a sequence $(r_k)_{k\in \mathbb{N}}\subset \mathbb{Q}^n$ such that
\[
\lim_{k\to\infty} r_k =r \quad \text{and} \quad \forall t>0:\  \lim_{k\to\infty} \sup_{s\in[-t,t]} \norm{g_k(r_k s)-h(s)}=0.
\] 
We call the sequence $(g_k)_k$ the \it{$\mathbb{Z}^n$-periodic sequence} of the function $h$.
\end{definition} 
\begin{remark}
Remark that for every $k\in \mathbb{N}$ the function $s\mapsto g_k(r_k s)$ is a  periodic function.
\end{remark}
 \begin{mainResultI}
Let   $f:\mathbb{R}^n\to\mathbb{R}^n$ be a trigonometric polynomial function. For every   $x_0\in\mathbb{R}^n$  the unique solution  $x:\mathbb{R}\to\mathbb{R}^n$     of the   differential  equation
\begin{equation}\label{edo:general}
\dot{x}=  f(x),\quad t\in\mathbb{R}, \quad  x(0)=x_0,
\end{equation}
admits  a rotation vector $\rho\in \mathbb{R}^n$. In addition, the function 
\[
t\mapsto x(t)-\rho t,
\]
is weakly almost periodic of slope $\rho$.
\end{mainResultI}
\section{Space of $C^\infty$    periodic modulo $\mathbb{Z}^n$ functions}
We define in this Section the space and the norm used   to prove the  Main result. The proofs of Lemmas for this Section are left in Appendix. In order to use the Fourier development,  let us introduce the following notation.
\begin{notation}\label{nota:ap}
For every  continuous function $g:\mathbb{R}^n\to\mathbb{R}^n$ and every $p\in 2\pi \mathbb{Z}^n$ we denote $ a_p[g]\in \mathbb{C}^n$  the following limit if it exists
\[
 a_p[g]:=\lim_{t\to+\infty}\frac{1}{(2t)^n}\int_{-t}^t  \ldots \int_{-t}^t  g(z)\exp( - i \langle z,p\rangle)dz_1\ldots dz_n.
\]
\end{notation}
  In this Section and the Section \eqref{mainprop}, for every function $g:\mathbb{R}^n\to\mathbb{R}^n$ and every $\alpha \in \mathbb{R}$,  we denote  $g_\alpha$ the function defined as $g_\alpha(z):=g(\alpha z)$ for all $z\in \mathbb{R}^n$.
 The following constant $\omega$ will be used    as change of variable in order to find a contraction in Lemmas \ref{lem:invariantcompact} and \ref{lem:beforprop} of Section \ref{mainprop}.   For every  $\omega\in\mathbb{N}^*$, we denote $E_\omega(\mathbb{R}^n)$ the set of $C^\infty$ function $g:\mathbb{R}^n \to \mathbb{R}^n$ such that $g_\omega$ is a  periodic modulo $\mathbb{Z}^n$ function. We remark, for very  $\omega\in\mathbb{N}^*$ and   $g\in E_\omega(\mathbb{R}^n)$  that
\begin{equation}\label{lem:transfourier1}
 a_p[g_\omega ]= \int_{0}^1 \ldots \int_{0}^1  g_\omega (  z)\exp( - i \langle z,p\rangle)dz_1\ldots dz_n,
\end{equation}
which is the Fourier coefficient of the function $g_\omega$.  Since $g\in E_\omega(\mathbb{R}^n)$  is $C^\infty$, by Dirichlet Theorem,  
\begin{gather}
 \label{lem:transfourier}  g (  z)=\sum_{p \in  2\pi \mathbb{Z}^n }a_{p }[g_\omega ] \exp(i\frac{1}{\omega} \langle z,p \rangle ),\\
\notag\forall q\geq 0:\quad \sum_{p \in 2\pi \mathbb{Z}^n}\norm{a_{p }[g_\omega] }\norm{\frac{p}{\omega}}^q<+\infty.
\end{gather} 
 We are now in position to define the following seminorm in $E_\omega(\mathbb{R}^n)$: Let   $\omega\in\mathbb{N}^*$ and    $g\in E_\omega(\mathbb{R}^n)$, we denote for every $\omega\in \mathbb{N}^*$ and $q\geq 0$
\[
\norm{   g }_{\omega,q}:=2\sum_{p \in 2\pi \mathbb{Z}^n/\mathbb{0}}\norm{a_{p }[g_\omega] }\norm{\frac{p}{\omega}}^q,
\]
where $\mathbb{0}:=(0,\ldots, 0)\in\mathbb{R}^n$ and where we recall that  $\norm{y}:=\max_{1\le j\le n}\norm{y_j}$ for every $y:=(y_j)_{j=1}^n\in\mathbb{C}^n$. 
We prove in the following Lemma that a periodic modulo $\mathbb{Z}^n$ function $g$ is $C^\infty$   if it is {\it  uniformly bounded for the seminorm}, i.e
\[
\forall q\geq 0:\quad \norm{    g }_{\omega,q}<+\infty,
\]
In other words, the set $E_\omega(\mathbb{R}^n)$  is include in  the set of the   periodic modulo $\mathbb{Z}^n$ functions  uniformly bounded for the seminorm.  
\begin{lemma}\label{lem:cinfty}
Let   $\omega\in\mathbb{N}^*$.   Let  $(c_p)_p$ be a complex-valued family such that
\[
\forall q\geq 0:\quad  \sum_{p \in 2\pi \mathbb{Z}^n/\mathbb{0}}\norm{c_{p }  }\norm{\frac{p}{\omega}}^q<+\infty.
\]
Then the following series is normally convergent
\[
g(z):=\sum_{p \in 2\pi \mathbb{Z}^n}c_p\exp( i \frac{1}{\omega}\langle z,p\rangle ), \quad c_p \in \mathbb{C}^n,
\]
and $c_p=a_p[g_\omega]$ for every $p\in 2\pi\mathbb{Z}^n$. Further,   $g\in E_\omega(\mathbb{R}^n)$.  
\end{lemma}
\begin{proof}
Appendix. A
\end{proof}
 In the following Lemma we prove that the seminorm $\norm{   .}_{\omega,0}$  is a norm on the space $\{g\in E_\omega(\mathbb{R}^n):\ g(\mathbb{0 })=\mathbb{0}\}$ and we compare it to the uniform norm topology.
\begin{lemma}\label{lem:equtow}
Let be  $\omega\in \mathbb{N}^*$ and $g\in E_\omega(\mathbb{R}^n)$     such that $g(\mathbb{0})=\mathbb{0}$ then
\[
\norm{a_{\mathbb{0}}[g_\omega]}\le \frac{1}{2}\norm{   g}_{\omega,0}\quad\text{and}\quad \norm{   g}_\infty \le  \norm{   g}_{\omega,0}.
\]
\end{lemma} 
\begin{proof}
Appendix. B
\end{proof}
We denote $d^kg$ the $k^{th}$ differential of a function $g:\mathbb{R}^n\to\mathbb{R}^n$. The following Lemma   gives an upper-bound of the quantity $\norm{d^kg}_{\omega,q}$ when $g$ is a trigonometric polynomial.   We recall that   $g_{\frac{1}{\omega}}(z):=g(\frac{z}{\omega})$.
\begin{lemma}\label{lem:beta}
Let $g :\mathbb{R}^n\to\mathbb{R}^n$ be a trigonometric polynomial function. Then there exists $\beta:=\beta(g)>0$ such that for every $\omega\in\mathbb{N}^*$ we have
\[
\norm{ d^k g_{\frac{1}{\omega}}}_{\omega,q}< n^k \beta (\frac{\beta}{\omega})^{q+k},\quad \forall q,k\geq0.
\]
\end{lemma}
\begin{proof}
Appendix. C
\end{proof}
We end this Section by the following inequality.
\begin{lemma}\label{lem:ineq}
Let   $\omega \in \mathbb{N}^*$ and  $\Big(h_j\in E_\omega(\mathbb{R}^n)\Big)_{j=1 }^k$ .  Then
\[
\forall q\in\mathbb{N},\ \forall k\in\mathbb{N}^*:\  \norm{   \Pi_{j=1}^k h_j }_{\omega,q}\le \frac{(k \omega^{k-1})^q}{2^{ k-1 }}   \Pi_{j=1}^k   [\norm{2 a_{\mathbb{0}}[h_{j,\omega}]}+ \norm{  h_j}_{\omega,q}].
\]
\end{lemma} 
\begin{proof}
Appendix. D
\end{proof}
\section{Main proposition}\label{mainprop}
The Main result affirms that the solution $x$ of Equation \eqref{edo:general} is a sum of a linear part and a bounded part. The strategy to prove the Main result  is to approximate the bounded part of  $x$ by  a  $\mathbb{Z}^n$-periodic sequence.   Using the Fourier development and  Equation \eqref{lem:transfourier},   remark,  for every $C^\infty$ periodic modulo $\mathbb{Z}^n$ function $g$  that
 \[
 f(z+g(z))= \sum_{p\in 2\pi \mathbb{Z}^n}  a_p[H[g]]  \exp( i  \langle z,p\rangle ),\ H[g](z):= f(z+g(z)),
 \]
under some convergence  assumption of  the series, by integration we get
\begin{align*}
\forall v\in \mathbb{R}^n:\quad \int_0^t &f(v s+g(v s))ds =t  \sum_{p\in 2\pi \mathbb{Z}^n,\ \langle v,p\rangle= 0 } a_p[H[g] ]  \\
&+\sum_{p\in 2\pi \mathbb{Z}^n,\ \langle v,p\rangle\neq 0 } \frac{ a_p[H[g]]] }{i \langle v,p\rangle } \Big(\exp( i   \langle v,p\rangle t )-1\Big).
\end{align*} 
The last term of the right member of the last equality  will play the role of  $\mathbb{Z}^n$-periodic sequence of the bounded part of the solution $x$ of Equation \eqref{edo:general}. In order to find an  upper-bound of the bounded part,  let us introduce the following notations.
\begin{notation}\label{Not:Tau} 
Let $f: \mathbb{R}^n \to \mathbb{R}^n$ be  a trigonometric polynomial. We denote the finite subset $\Lambda_f\subset 2\pi\mathbb{Z}^n$ as
\[
\Lambda_f:=\{p\in 2\pi\mathbb{Z}^n:\quad \norm{a_p[f]}\neq0\},
\]
and we denote 
\[
|\Lambda_f|:=\max\{\norm{p},\quad p\in \Lambda_f \}.
\]
Let be $y\in  \mathbb{R}^n/\{\mathbb{0}\} $. Define,
\[
\Lambda(f,y):=  \{ p\in   2\pi\mathbb{Z}^n:\    \norm{p} \le 2\pi+|\Lambda_f|, \ \langle y ,p\rangle\neq0\}, 
\]
Remark that  $\Lambda(f,y)\neq \emptyset$. We denote
\[ 
\tau( f, y) :=  \max \Big\{\frac{1}{ |\langle y ,p\rangle |}:\quad\ p\in \Lambda(f,y)\}.
\]
Let be $y\in\mathbb{Q}^n/\{\mathbb{0}\}$, we denote
\[
\tau(  y)=\max\Big\{\frac{1}{ |\langle y ,p\rangle | }:\quad p\in  2\pi\mathbb{Z}^n,\ \langle y ,p\rangle\neq0\Big\}.
\] 
We denote $\beta$ the constant $\beta(f)$ of the  function $f$ defined in Lemma  \ref{lem:beta}.
\end{notation}
In the following Proposition we prove that the bounded part of the solution $x$ of Equation \eqref{edo:general} can be approximated by a $\mathbb{Z}^n$-periodic functions and we find an appropriate upper-bound.
\begin{proposition}\label{lem:doublefixedpoint}[Main proposition]
Let $f :\mathbb{R}^n\to\mathbb{R}^n$ be a trigonometric polynomial function.  Then for every $r\in  \mathbb{Q}^n/\{\mathbb{0}\}$ and every $\epsilon>0 $ there exists a $C^\infty$ periodic modulo $\mathbb{Z}^n$ function $\phi_{r,\epsilon}:\mathbb{R}^n\to\mathbb{R}^n$     such that  $\norm{\phi_{r,\epsilon}}_\infty  < 2\beta  \tau(f,r )$  and such that
\begin{equation*}
\sup_{z\in \mathbb{R}^n}\norm{\phi_{r,\epsilon}(z)- \sum_{p\in 2\pi\mathbb{Z}^n ,\ \langle r,p\rangle\neq0} \frac{ a_p[H[\phi_{r,\epsilon}]]}{ i\langle r,p\rangle} (\exp(i \langle z,p\rangle )-1)}<\epsilon,
\end{equation*}
where
\[
  H[\phi_{r,\epsilon}](z):=f( z+  \phi_{r,\epsilon}(z)),\ \forall z\in\mathbb{R}^n.
\]  
\end{proposition}  
As is state in the above Section, the following constant $\omega$ is used as change of variable in order to find a contraction.
\begin{definition}
For  every $r\in  \mathbb{R}^n$ and every   $\omega\in \mathbb{N}^*$, define the  set  $K_{ r,\omega}$   as $g\in K_{  r,\omega}$ if 
\begin{itemize} 
\item  there exists a complex-valued family $(c_p)_{p\in 2\pi\mathbb{Z}^n}$ such that
\[
g(z)= \sum_{p\in   2\pi \mathbb{Z}^n } c_p ( \exp(i \frac{1}{\omega} \langle z,p\rangle)-1),\ \forall z\in \mathbb{R}^n,
\]
\item $\norm{ g}_{\omega,0} \le2 \beta  \tau(f,r )$,
\item $\norm{ g}_{\omega,q}<\infty$ for every $q\geq1$.
\end{itemize}
\end{definition}
\begin{lemma}\label{lem:kine}
The set $K_{ r,\omega}$ is a nonempty subset of   $E_\omega(\mathbb{R}^n)$. 
\end{lemma}
\begin{proof}
The set $K_{ r,\omega}\neq \emptyset$ because it contains the function $z\mapsto g(z)=\mathbb{0}$.   By definition of $K_{ r,\omega}$ and by Lemma \ref{lem:cinfty} the function $g$ is $C^\infty$.
\end{proof}
For  every $r\in  \mathbb{R}^n$ , for every   $\omega\in\mathbb{N}^*$, and every $g\in K_{ r,\omega}$ let be  ${\Psi} [  r,\omega ,g] $ the function defined by the following series in its convergence  domain
\begin{align*} 
\forall z\in \mathbb{R}^n:\quad \Psi [ r,\omega ,g](   z)&  :=    \sum_{p\in   2\pi \mathbb{Z}^n,\  \langle r,p\rangle\neq0} \frac{ a_p[H[\omega,g] ]}{ i\langle r,p\rangle} (\exp(i \frac{1}{\omega} \langle z,p\rangle)-1),
\end{align*} 
where 
\[
H[\omega,g](z):=f\Big(  z+  g(\omega z) \Big) .
\] 
Since $f$ is a real polynomial trigonometric function then ${\Psi}[  r,\omega ,g](z)\in \mathbb{R}^n$ for every $z\in \mathbb{R}^n$ such that the  series converge.
In the following Lemma we prove that $K_{  r,\omega}$ is invariant under the operator ${\Psi} [ r,\omega ,.]$. We deduce that ${\Psi} [ r,\omega ,g](z)$ is defined for every $g\in K_{  r,\omega}$ and $z\in \mathbb{R}^n$.
\begin{lemma}\label{lem:invariantcompact}
Let $f :\mathbb{R}^n\to\mathbb{R}^n$ be a trigonometric polynomial function.  For  every $r\in    \mathbb{Q}^n/\{\mathbb{0}\}$, there exists $\omega_r>0$ such that for every $\omega>\omega_r$ we have
\[
 g\in K_{  r,\omega}\implies  {\Psi} [ r,\omega ,g]\in K_{  r,\omega}.
 \]
In addition, ${\Psi} [ r,\omega ,g]$ is defined for every $z\in \mathbb{R}^n$.
\end{lemma}
\begin{proof}
Prove that   
\[
\forall g\in K_{  r,\omega}:\quad  \norm{{\Psi} [ r,\omega ,g]}_{\omega,q}<\infty ,\quad \forall q\geq0.
\]
 Let be $g\in K_{  r,\omega}$   and denote 
\[
 \tilde{H}[\omega,g](z):= f\Big(  z+ g(\omega z) \Big)-f (z ).
\] 
We have
\[
a_p[H[\omega,g]]=a_p[f  ]+a_p[\tilde{H}[\omega,g]].
\]
By definition of ${\Psi} [ r,\omega ,.]$ for every $g\in K_{  r,\omega}$,
\begin{align*} 
\Psi [ r,\omega ,g](   z)&   =   \sum_{p\in   2\pi \mathbb{Z}^n,\  \langle r,p\rangle\neq0} \frac{ a_p[f  ]}{ i\langle r,p\rangle} (\exp(i \frac{1}{\omega} \langle z,p\rangle)-1)\\
&+   \sum_{p\in    2\pi \mathbb{Z}^n,\  \langle r,p\rangle\neq0} \frac{ a_p[\tilde{H}[\omega,g]]}{ i\langle r,p\rangle} (\exp( i\frac{1}{\omega} \langle z,p\rangle)-1)
\end{align*}
Recall that $|\Lambda_f|$ is defined on the Notations \ref{Not:Tau}:  For every $p\in   2\pi \mathbb{Z}^n$ such that $ \langle r,p\rangle\neq0$ and such that $\norm{p}> |\Lambda_f|$ we have  $a_p[ f ]=0$.  By definition of $ \tau(f,  r)$ and $ \tau(  r)$ in Notation \ref{Not:Tau}, we get  
\begin{align}
\notag \norm{ \Psi [r,\omega,g]}_{\omega,q}&\le      2\tau( f, r)  \sum_{p\in 2\pi \mathbb{Z}^n,\  \langle r,p\rangle\neq0 }\norm{a_p [ f ]}\norm{\frac{p}{\omega}}^q  \\
\notag&+  2\tau( r)    \sum_{p\in 2\pi \mathbb{Z}^n ,\  \langle r,p\rangle\neq0}\norm{a_p[ \tilde{H}[\omega,g]] }   \norm{\frac{p}{\omega}}^q.
\end{align} 
By definition of the seminorm, 
\begin{align}
\label{equ:dpsi} \norm{ \Psi [r,\omega,g]}_{\omega,q} \le    \tau( f, r)   \norm{ f_{\frac{1}{\omega }}}_{\omega,q}+  \tau( r)    \norm{  \Big(\tilde{H}[\omega,g]\Big)_{\frac{1}{\omega}}}_{\omega,q}.
\end{align} 
By Lemma \ref{lem:beta}, we have
\begin{align}
\label{equ:dpsi1}  \norm{ f_{\frac{1}{\omega }} }_{\omega,q} \le    \frac{\beta^{ q+1}}{\omega^{ q}},\quad \forall q\geq0 .
\end{align}
Now,  estimate the quantity    $ \norm{  \Big(\tilde{H}[\omega,g]\Big)_{\frac{1}{\omega}}}_{\omega,q}$.  By definition,
\[
 \Big(\tilde{H}[\omega,g]\Big)_{\frac{1}{\omega}}(z)= f\Big(\frac{1}{\omega}  z+ g(z )\Big)-f (\frac{1}{\omega} z )=\sum_{k=1}^\infty \frac{d^kf(\frac{z}{\omega} )}{k!} (g(z))^{(k)},
\]
where 
\[
 d^kf(\frac{z}{\omega}) (g( {z} ))^{(k)} :=\sum_{i_1,\ldots,i_k=1}^n\frac{\partial f  }{\partial z_{i_1}\ldots \partial z_{i_k}}(\frac{z}{\omega})    g_{i_1}(z)\ldots g_{i_k}(z).
\] 
Since $f$ is polynomial trigonometric function, then
\[
\forall z\in\mathbb{R}^n:\quad f_{\frac{1}{\omega}}(z)=f(\frac{z}{\omega}) =a_{\mathbb{0}}[f]+\sum_{p\in 2\pi\mathbb{Z}^n, p\neq \mathbb{0}  }a_p[f]\exp( i\frac{1}{\omega}\langle z,p\rangle ),
\]
then for every $k\geq1$ we get
\begin{align*}
d^{k} f_{\frac{1}{\omega}}(z)&=\sum_{p\in 2\pi\mathbb{Z}^n, p\neq \mathbb{0}  }a_p[f]d^{k} \Big[\exp( i\frac{1}{\omega}\langle z,p\rangle )\Big]\\
&=\sum_{p\in 2\pi\mathbb{Z}^n, p\neq \mathbb{0}  }\Big(a_p[f]\sum_{s_1=1}^n\ldots \sum_{s_k=1}^n i^k \frac{p_{s_1}}{\omega}\ldots \frac{p_{s_k}}{\omega}\Big)\exp( i \frac{1}{\omega}\langle z,p\rangle),
\end{align*}
By Notation \ref{nota:ap}
\[
 a_p[d^{k} f_{\frac{1}{\omega}}]:=\lim_{t\to+\infty}\frac{1}{(2t)^n}   \sum_{p\in 2\pi\mathbb{Z}^n, p\neq \mathbb{0}  }\Big(a_p[f]\sum_{s_1=1}^n\ldots \sum_{s_k=1}^n i^k \frac{p_{s_1}}{\omega}\ldots \frac{p_{s_k}}{\omega}\Big)\theta(\omega,-t,t),
\]
where
\[
\forall t\in \mathbb{R}:\ \theta(\omega,t,-t):= \int_{-t}^t  \ldots \int_{-t}^t  \exp( i \frac{1}{\omega}\langle z,p\rangle)dz_1\ldots dz_n.
\]
Since  $p\neq \mathbb{0}$, then
\[
\lim_{t\to+\infty}\frac{1}{(2t)^n} \theta(\omega,t,-t)=0,
\]
we deduce that 
\[
 \norm{   a_\mathbb{0} [d^{k} f_{\frac{1}{\omega}}  ]}=0,\quad \forall k\geq1.
\]
Since  $g\in K_{r,\omega}$ then $g(\mathbb{0})=\mathbb{0}$. Thanks to Lemma \ref{lem:equtow} we obtain
\[
\norm{ a_\mathbb{0} [g ]}    \le \norm{  g }_{\omega,0}\le \alpha_0:=2 \beta  \tau(f,r ).
\] 
By Lemma \ref{lem:ineq}, we have for all $q\geq0$
\begin{align}
\notag   \norm{  \Big(\tilde{H}[\omega,g]\Big)_{\frac{1}{\omega}}}_{\omega,q} &\le \sum_{k=1}^\infty \frac{1}{k!}    \norm{ (d^k f_{\frac{1}{\omega}})(g )^{(k)} }_{\omega,q}\\
\notag &\le \sum_{k=1}^\infty \frac{(k+1)^q \omega^{kq} }{k!2^k}    \norm{   d^{k} f_{\frac{1}{\omega}} }_{\omega,q}(\norm{  g}_{\omega,q}+\alpha_0)^{  k}.
\end{align} 
By Lemma \ref{lem:beta}, we find
\begin{align*}
  \forall q\geq0:\quad   \norm{  \Big(\tilde{H}[\omega,g]\Big)_{\frac{1}{\omega}}}_{\omega,q} &\le  \sum_{k=1}^\infty \frac{\exp(qk) \omega^{kq} }{k!2^k}     n^k  \frac{\beta^{k+q+1}}{\omega^{k+q}}   [\alpha_0+\norm{  g }_{\omega,q}]^{  k}.
\end{align*}
By hypothesis $g\in K_{  r,\omega}$, then
\begin{equation*}
\forall q\geq0, \exists \alpha_q>0:\quad \norm{     g}_{\omega,q} \le \alpha_q \ \text{where}\  \alpha_0:=2 \beta \tau( f, r).
\end{equation*} 
 we deduce that for all $q\geq0$, 
\begin{align*}
\notag \norm{  \Big(\tilde{H}[\omega,g]\Big)_{\frac{1}{\omega}}}_{\omega,q}  &\le  \sum_{k=1}^\infty \frac{\exp(qk) \omega^{kq}  }{k!2^k}    n^k    \frac{\beta^{k+q+1}}{\omega^{k+q}}  (\alpha_0+\alpha_q)^{  k}\\
&= \frac{\beta^{ q+1}}{\omega^{ q}} \sum_{k=1}^\infty \frac{1 }{k!}    \Big(n \exp(q ) \omega^{ q} \frac{\beta  }{2\omega }  (\alpha_0+\alpha_q)\Big)^{  k}\\
&\le \frac{\beta^{ q+1}}{\omega^{ q}} \Big[ \exp\Big(n  \exp(q ) \omega^{ q} \frac{\beta  }{2\omega }  (\alpha_0+\alpha_q)\Big)-1\Big]<\infty.
\end{align*}
By Equations   \eqref{equ:dpsi1} and  \eqref{equ:dpsi}, we  obtain
\[
\forall g\in K_{  r,\omega}:\quad  \norm{{\Psi} [ r,\omega ,g]}_{\omega,q}<\infty ,\quad \forall q\geq0.
\]
Choose $ \omega>\omega_r>0$, where $\omega_r\in \mathbb{N}^*$ satisfies 
\[
 \beta\Big[\exp\Big(n   \frac{\beta  }{\omega_r }  \alpha_0 \Big)-1\Big]< \frac{\beta       \tau(f,r )}{ \tau(r)}.
\]
We obtain 
 \begin{align}
\label{equ:dpsi2}\forall \omega >\omega_r:\quad  \norm{   \Big(\tilde{H}[\omega,g]\Big)_{\frac{1}{\omega}}}_{\omega,0}  <\frac{\beta       \tau(f,r )}{ \tau(r)}.
\end{align} 
Replace both Equations \eqref{equ:dpsi1} and \eqref{equ:dpsi2} on Equation \eqref{equ:dpsi}, we obtain
\begin{align*}
\forall \omega >\omega_r:\quad \norm{  \Psi [r,\omega,g]}_{\omega,0} &<  \beta  \tau(f, r)  +  \beta \tau(f, r)    = 2\beta \tau( f, r)=\alpha_0.
\end{align*}
\end{proof}

\begin{lemma}\label{lem:beforprop}
Let $f :\mathbb{R}^n\to\mathbb{R}^n$ be a trigonometric polynomial function. For  every $r\in    \mathbb{Q}^n/\{\mathbb{0}\}$ there exists $\omega:=\omega(r) >0$ such that for every $\epsilon>0$, there exists    $\phi_{r,\omega, \epsilon}  \in  K_{ r,\omega} $   satisfying
\[
\norm{\phi_{r,\omega, \epsilon}-{\Psi} [r,\omega  ,\phi_{r,\omega, \epsilon}]}_\infty<\epsilon.
\]
\end{lemma}
\begin{proof}
 Let be $\omega>\beta$ and $h,g \in K_{ r,\omega}$.  For every fixed $s\in [0,1 ]$ , define the function $V_s \in K_{ r,\omega}$ as
\[
V_s(z):=  s h(z)+(1-s) g(z ) ,\quad \forall z\in \mathbb{Z}^n.
\] 
\begin{align*}
\forall z\in \mathbb{R}^n:\quad \Psi[r,\omega,h](z)-\Psi[r,\omega,g](z)&=\Psi[r,\omega,V_1(z)]-\Psi[r,\omega,V_0(z)]\\
&=\int_0^1 \frac{d}{ds}\Psi[r,\omega,V_s(z)]ds.
\end{align*}
By definition of $\Psi$ we get
\begin{gather*}
  \Psi[r,\omega,h](z)-\Psi[r,\omega,g](z) \\
=\int_0^1  \sum_{p\in   2\pi \mathbb{Z}^n,\  \langle r,p\rangle\neq0}\frac{ \frac{d}{ds}a_p[H[\omega,V_s]] }{ i\langle r,p\rangle} (\exp(i \frac{1}{\omega} \langle z,p\rangle)-1)ds.
\end{gather*}
We have 
\[
\frac{d}{ds}a_p[H[\omega,V_s]]  = a_p[  \frac{d}{ds} H[\omega,V_s]],\quad H[\omega,V_s](z):=f\Big(  z+  V_s(\omega z) \Big).
\] 
Since $\frac{d}{ds}V_s=h-g$, then
\begin{align*}
 \frac{d}{ds} H[\omega,V_s](z) &=df\Big( z+V_s(\omega z) \Big)  \frac{d}{ds}V_s(\omega z)\\
&=df\Big(  z+V_s (\omega z) \Big)[h-g](\omega z).
\end{align*}
Then
\[
a_p[  \frac{d}{ds}H[\omega,V_s]]  =  a_p\Big[ \phi[\omega,V_s] (h-g)_\omega \Big],
\]
where
\begin{equation}\label{equ:phiomegavs}
\phi[\omega,V_s](z):=   df\Big(  z+V_s(\omega z) \Big) .
\end{equation}
For every fixed $s\in [0,1 ]$ we have  $\phi[\omega,V_s]\in K_{ r,\omega}$. By  Lemma \ref{lem:ineq}, for every fixed $s\in [0,1]$ we have
\begin{align*}
  \norm{ \Big(\phi[\omega,V_s] (h-g)_\omega\Big)_{\frac{1}{\omega}}}_{\omega,0} \le    \norm{\Big(\phi[\omega,V_s]\Big)_{\frac{1}{\omega}} }_{\omega,0} \norm{ h-g}_{\omega,0}.
\end{align*}
Then  
\begin{align}
\label{equ:pointfixpsi}&\norm{ \Psi[r,\omega,h] -\Psi[r,\omega,g] }_{\omega,0} \le   \tau(  r) \Big(\sup_{s\in [0,1]}  \norm{\Big(\phi[\omega,V_s]\Big)_{\frac{1}{\omega}}}_{\omega,0}  \Big)\norm{  h-g}_{\omega,0} .
\end{align}
Prove that there exists $\omega_r>0$ such that for every $\omega>\omega_r$ we have
\begin{equation}\label{equa:half}
\forall g,h\in K_{r,\omega}:\quad  \tau(  r)\Big( \sup_{s\in[0,1]} \norm{\Big(\phi[\omega,V_s]\Big)_{\frac{1}{\omega}}}_{\omega,0}\Big) <\frac{1}{2}.
\end{equation}
As in  Proof of Lemma \ref{lem:invariantcompact}: By Equation \eqref{equ:phiomegavs}   we have for every$s\in [0,1]$,
\begin{align*}
\Big(\phi[\omega,V_s]\Big)_{\frac{1}{\omega}}&=  \sum_{k=0}^\infty \frac{ d^{k+1}f_{\frac{1}{\omega} }(z)}{k!} \Big( V_s(s)\Big)^{(k)} ,
\end{align*}
By Lemmas \ref{lem:ineq};
\begin{align*}
\forall s\in[0,1]:\ \norm{\Big(\phi[\omega,V_s]\Big)_{\frac{1}{\omega}}}_{\omega,0} \le    \sum_{k=0}^\infty  \frac{1}{k!2^k}   \norm{  d^{k+1}f_{\frac{1}{\omega}} }_{\omega,0}\norm{ V_s }_{\omega,0}^{2k}.
\end{align*}
By hypothesis $g,h \in K_{  r,\omega}$, by consequence
\[
 \sup_{s\in[0,1]}\norm{  V_s}_{\omega,0}\le \norm{ g}_{\omega,0}+\norm{  h}_{\omega,0}\le 4\beta \tau(f,r ).
\]
Using Lemma\ref{lem:beta}, we get
\begin{align}
\notag \sup_{s\in[0,1]}\norm{ \Big(\phi[\omega,V_s]\Big)_{\frac{1}{\omega}}}_{\omega,0}&\le \sum_{k=0}^\infty \frac{1}{k!2^k}  \frac{n^{k+1} \beta^{k+2}}{\omega^{k+1}}(4\beta \tau(f,r ))^k \\
\notag &=  \frac{n\beta^2}{\omega }\sum_{k=0}^\infty \frac{1}{k!} \Big(\frac{2n \beta^2 }{\omega }  \beta \tau(f,r )\Big) ^{ k}\\
\notag &= \frac{n\beta^2}{\omega } \exp\Big(\frac{2n \beta^2 }{\omega }  \beta \tau(f,r )\Big).
\end{align}
Choose $\omega_r>0 $ large such that 
\[
 \tau(  r)    \frac{n\beta^2}{\omega_r } \exp\Big(\frac{2n \beta^2 }{\omega_r }  \beta \tau(f,r )\Big)<\frac{1}{2}.
\] 
We have proved Equation \eqref{equa:half}. Thanks to Equation \eqref{equ:pointfixpsi},   for every $\omega>\omega_r$   we have
\begin{equation*}
  \norm{ [\Psi[r,\omega,h] -\Psi[r,\omega,g]}_{\omega,0}\le \frac{1}{2}  \norm{ h-g}_{\omega,0}.
\end{equation*}
Now, choose $\omega >0$ fixed and   large.  Let be $g\in K_{r,\omega}$,  by  the last inequality,  for every $\epsilon>0$ there exists $k_\epsilon \geq0$  such that
\begin{align*}
 \norm{\Psi^{k_\epsilon+1}[ r,\omega , g]-\Psi^{k_\epsilon}[ r,\omega , g]}_{\omega,0} <\epsilon,
\end{align*}
Denote
\[
\phi_{r,\omega, \epsilon}:=\Psi^{k_\epsilon}[ r,\omega , g].
\]
By Lemma \ref{lem:invariantcompact} we have $\phi_{r,\omega, \epsilon}\in K_{r,\omega}$. By Lemma \ref{lem:equtow} we obtain 
\begin{align*}
 \norm{\Psi[ r,\omega ,\phi_{r,\omega, \epsilon}]-\phi_{r,\omega, \epsilon}}_\infty <\epsilon,
\end{align*}

\end{proof}
\begin{proof}[Proof of Proposition \ref{lem:doublefixedpoint}]
Let be $r\in   \mathbb{Q}^n/\{\mathbb{0}\}$. By Lemma \ref{lem:beforprop}, there exists $\omega:=\omega_r >0$  such that for every $\epsilon >0$ there exists   $\phi_{r,\omega, \epsilon}\in K_{ r,\omega}$  satisfying
\begin{equation}\label{equ:proofpropfix}
 \norm{\Psi[  r,\omega, \phi_{r,\omega, \epsilon}]-\phi_{r,\omega, \epsilon}}_\infty<\epsilon.
\end{equation}
Define the functions,
\[
\tilde{\phi}_{r, \epsilon}(z):=\  \phi_{r,\omega, \epsilon}(\omega z),\ \text{and}\  H[g](z):=f( z+  g(z)).
\]
We recall that, 
\[
 H[\omega, \phi_{r,\omega, \epsilon}](z) =f\Big(  z+  \phi_{r,\omega, \epsilon}(\omega z) \Big).
\] 
By Equation \eqref{lem:transfourier1}, 
\[
 a_p[H[\omega,\phi_{r,\omega, \epsilon}]]= a_p[H[\tilde{\phi}_{r, \epsilon}]], 
\]
Using the definition of $\Psi[  r,\omega, \tilde{\phi}_{r, \epsilon}]$ and replace on Equation \eqref{equ:proofpropfix}, the function $\tilde{\phi}_{r, \epsilon}$ satisfies
\[
\sup_{z\in \mathbb{R}^n}\norm{\tilde{\phi}_{r, \epsilon}(z)- \sum_{p\in i2\pi\mathbb{Z}^n ,\ \langle r,p\rangle\neq0} \frac{ a_p[H[\tilde{\phi}_{r, \epsilon}]]}{ i\langle r,p\rangle} (\exp( i\langle z,p\rangle )-1)}<\epsilon.
\] 
By Lemma \ref{lem:beforprop}, $\phi_{r,\omega, \epsilon}\in K_{ r,\omega}$  then
\[
   \norm{\tilde{\phi}_{r, \epsilon}]}_\infty\le    \norm{\phi_{r,\omega, \epsilon}}_\infty <  2 \beta  \tau(f,r ).
\]
By Lemma \ref{lem:kine} the set $K_{ r,\omega}$ is a subset of $E_\omega (\mathbb{R}^n)$.  Then $\tilde{\phi}_{r, \epsilon}$ is a $C^\infty$ periodic modulo $\mathbb{Z}^n$ function.

\end{proof}

\section{Proof of the Main result}
 
\begin{proof}[Proof of Main results]
Consider the System \eqref{edo:general} where $f$ is a polynomial function.There exists $q\in \mathbb{N}^*$ such that $\norm{f}_\infty<q$. Use the change of variables
\[
x_q(t)=x(t)-x_0+qt,\quad \forall t\in \mathbb{R}, 
\]
we get
\begin{equation}\label{edo:generalB}
\frac{d}{dt}{x}_q(t)=  f(x_q(t)+x_0-q t\mathbb{1})+q,\quad t\in\mathbb{R}, \quad  x_q(0)=\mathbb{0}.
\end{equation}
where   $\mathbb{1}:=(1,\ldots, 1)\in\mathbb{R}^n$. Now, transform the last system to an autonomous systems. Define the functions $x_{n+1}:\mathbb{R}\to \mathbb{R}$ as the identity function: $x_{n+1}(t):=t$ for every $t\in \mathbb{R}$, the system \eqref{edo:generalB} can be written as
\begin{align*} 
\dot{x}_q&=  f(x_q+x_0-q x_{n+1}\mathbb{1})+q,\quad t\in\mathbb{R}, \quad x_q(0)=\mathbb{0},\\
\dot{x}_{n+1}&=1\quad t\in\mathbb{R}, \quad  x_{n+1}(0)=0,
\end{align*}
in other words,
\begin{align*} 
\dot{\tilde{x}} =  f_q(\tilde{x}),\quad t\in\mathbb{R}, \quad \tilde{x}=(x_q, x_{n+1}),\quad   \tilde{x}(0)=(\mathbb{0},0), 
\end{align*}
where $f_q:\mathbb{R}^{n+1}\to\mathbb{R}^{n+1}$  satisfies
\[
f_q(\tilde{z}):=\Big(f(z+x_0-q z_{n+1}\mathbb{1}),1\Big),\quad \forall \tilde{z}:=(z,z_{n+1})\in \mathbb{R}^{n+1}.
\]
Since $q\in \mathbb{N}^*$ then $f_q$ is a polynomial trigonometric function. In addition, $\min_{z\in \mathbb{R}^{n+1}}f_q(z)\geq \min\{q-\norm{f}_\infty, 1\}>0$. Without loss of generality, we consider the system \eqref{edo:general} by supposing that  $x(0)=\mathbb{0}$ and $f(z)\neq\mathbb{0}$ for all $z\in\mathbb{R}^n$.  Let  $(\epsilon_k)_k\subset (0,1]$ be a sequence satisfying $\lim_{k\to\infty}\epsilon_k=0$.  For every $k\geq 1$ let    $\gamma_k\in\mathbb{R}^n\to   \mathbb{Q}^n/\{\mathbb{0}\}$ be a   function satisfying-
\[
 \forall y\in \mathbb{R}^n:\quad \norm{y-\gamma_k(y)}<\epsilon_k. 
\]
We have $ \gamma_k(y) \in   \mathbb{Q}^n/\{\mathbb{0}\}$. For every $y\in\mathbb{R}^n$ and every $k\geq1$  consider  the function ${\phi}_{\gamma_k(y)}$ satisfying    the Main Proposition such that 
\begin{equation*}
\sup_{z\in \mathbb{R}^n}\norm{{\phi}_{\gamma_k(y)}(z)-\sum_{p\in 2\pi\mathbb{Z}^n ,\ \langle \gamma_k(y),p\rangle\neq0} \frac{ a_p[H[ \phi_{\gamma_k(y)}]]}{ i\langle \gamma_k(y),p\rangle} (\exp( i\langle  z,p\rangle )-1)}<\frac{1}{k},
\end{equation*}
Define the recurrent sequence $(\rho_k)_k\subset \mathbb{R}^n$ as
\[
\rho_0=\mathbb{0},\quad \rho_{k+1}:=\sum_{p\in  i2\pi\mathbb{Z}^n,\  \langle \gamma_k({\rho_k}),p\rangle =0  }   a_p[H[ {\phi}_{\gamma_k({\rho_k})}]],\quad \forall k\geq0.
\]
Prove that the sequence $(\rho_k)_k$ is bounded. Let be $\psi_k:\mathbb{R}\to\mathbb{R}^n$ the function defined by $t\mapsto \psi_k(t):={\phi}_{\gamma_k({\rho_k})}(\gamma_k({\rho_k})t)$. By     the Main Proposition, we get 
\begin{equation*}
\sup_{t\in \mathbb{R}}\norm{\psi_k(t)-\sum_{p\in  2\pi\mathbb{Z}^n ,\ \langle \gamma_k({\rho_k}),p\rangle\neq0} \frac{ a_p[H[ {\phi}_{\gamma_k({\rho_k})}]]}{ i \langle \gamma_k({\rho_k}),p\rangle} (\exp( i\langle \gamma_k({\rho_k}),p\rangle t)-1)}<\frac{1}{k},
\end{equation*} 
since the sum is normally convergent, that  implies
\[
\sup_{t\in \mathbb{R}}\norm{ {\psi}_k(t)- \int_0^t \sum_{p\in   2\pi\mathbb{Z}^n,\  \langle \gamma_k({\rho_k}),p\rangle \neq0  }  {a_p[H[ {\phi}_{\gamma_k({\rho_k})}]]} \exp(i \langle \gamma_k({\rho_k}),p\rangle s)ds}<\frac{1}{k}.
\]  
By Equation \eqref{lem:transfourier}, we have the following  Fourier development
 \[
 f(z+{\phi}_{\gamma_k({\rho_k})}(z))= \sum_{p\in 2\pi \mathbb{Z}^n}  a_p[H[{\phi}_{\gamma_k({\rho_k})}]]  \exp( i  \langle z,p\rangle ),
 \]
then
\begin{equation}\label{equ:fk}
\sup_{t\in \mathbb{R}}\norm{ {\psi}_k(t) -\int_0^t f(\gamma_k({\rho_k})s+ {\psi}_k (  s))ds-t \rho_{k+1}}<\frac{1}{k}.
\end{equation}
Since $\norm{\psi_k}_\infty<\infty$ then 
  \[
 \norm{\rho_{k+1} - \lim_{t\to\infty}\frac{1}{t}\int_0^t f(\gamma_k({\rho_k}) s+ {\psi}_k (  s))ds}=0.
 \] 
 we deduce that $\limsup_{k\to\infty}\norm{\rho_k}\le \norm{f}$. There exists  $\rho\in \mathbb{R}^n$ and a sub-sequence $(\rho_{k_s})_s$ which converge to $\rho$. In order to simplify the notation, we suppose that  $(\rho_{k})_k$   converge to $\rho$. 
  Since $\epsilon_k\to0$ then
\[
\lim_{k\to\infty}\rho_k=\lim_{k\to\infty}\gamma_k({\rho_{k}})=\rho.
\]
We have supposed  that $f(z)\neq \mathbb{0}$ for every $z\in\mathbb{R}^n$, then $\rho\neq \mathbb{0}$. There exists $c>0$ and $k_0\geq0$ such that
 \[
 \tau(f,{\gamma_k({\rho_k})})<c,\quad \forall k\geq k_0.
 \]
By   the Main Proposition, we obtain
\[
\sup_{k\geq k_0} \norm{ {\phi}_{\gamma_k({\rho_k})}}_\infty \ \le   2 \beta  \sup_{k\geq k_0} \tau(f,{\gamma_k({\rho_k})})<2 \beta c.
\]
 Now, prove that the sequence functions $(\psi_k)$ converge uniformly on every interval $[0,T]$. Since $f$ is a polynomial trigonometric function, then there exist $\eta>0$ such that $f$ is uniformly $\eta$-Lipschitz function.   For every $T>0$ we have
\begin{align*}
\sup_{t\in [0,T]}\norm{\psi_{k_2}(t)-\psi_{k_1}(t)}&=\sup_{t\in [0,T]}\norm{\exp( 2\eta t) \exp(-2\eta t) \psi_{k_2}(t)-\psi_{k_1}(t)}\\
&\le\exp( 2\eta T) \norm{\psi_{k_2} -\psi_{k_1}}_T,\quad \forall k_1,k_2\in \mathbb{N},
\end{align*}
where
\[
\norm{\psi_{k_2}-\psi_{k_1}}_T:=\sup_{t\in [0,T]}\norm{\exp(-2\eta t) [\psi_{k_2}(t)-\psi_{k_1}(t)]}.
\] 
It is sufficient to prove that
\[
 \lim_{k_2, k_1\to+\infty}\norm{\psi_{k_2} -\psi_{k_1 }}_T=0.
\]
By Equation \eqref{equ:fk}
\begin{align*}
 \norm{\psi_{k_2} -\psi_{k_1 }}_T & \le \eta \sup_{t\in [0,T]}\Big[ \exp(-2\eta t) \int_0^t s\norm{\gamma_{k_2}( \rho_{k_2 })-\gamma_{k_1}( \rho_{k_1 })}  ds\Big]\\
&+ \eta \sup_{t\in [0,T]} \Big[\exp(-2\eta t) \int_0^t \norm{ {\psi}_{k_2} (  s)-  {\psi}_{k_1} (  s)} ds\Big] \\
&+  \sup_{t\in [0,T]}\Big[ \exp(-2\eta t) t \norm{\rho_{k_2+1}- \rho_{k_1+1}}+\frac{1}{k_2}+\frac{1}{k_1}\Big]\\
&\le   \eta T^2 \norm{\gamma_{k_2}( \rho_{k_2 })-\gamma_{k_1}( \rho_{k_1 })} +\frac{1}{2} \norm{ {\psi}_{k_2} -  {\psi}_{k_1} }_T\\
&+T \norm{\rho_{k_2+1}- \rho_{k_1+1}}+\frac{1}{k_2}+\frac{1}{k_1}.
\end{align*}
Then
\begin{align*}
\frac{1}{2} \norm{ {\psi}_{k_2} -  {\psi}_{k_1} }_T&\le  \eta  T^2   \norm{\gamma_{k_2}( \rho_{k_2 })-\gamma_{k_1}( \rho_{k_1 })} \\
&+T\norm{\rho_{k_2+1}- \rho_{k_1+1}}+\frac{1}{k_2}+\frac{1}{k_1} \to 0, \ \text{when}\ k_2,k_1\to +\infty.
\end{align*}   
We deduce that the sequence function $(\psi_k)_k$ is a Picard iteration for  the solution of  the differential equation 
\[
\dot{x}=f(x),\quad x(0)=\mathbb{0},
\]
 there exists  a  weakly almost periodic function $\psi_{\rho}^*:\mathbb{R}\to\mathbb{R}^n$  of slope $\rho$  such that
\begin{align*}
\psi_{\rho}^* (t)&=\lim_{k\to\infty}\psi_k(t) = \int_0^t \lim_{k\to\infty} f({\gamma_k({\rho_k})}s+ \psi_k(t))ds-t \lim_{k\to\infty} \rho_k \\
&=\int_0^t   f(\rho s+\psi_{\rho}^* (  s))ds-t \rho,\ \forall t\in \mathbb{R}.
\end{align*}
By uniqueness of solution of differential equation, we have proved that
\[
x(t)= \rho t+ \psi_{\rho}^* (  t),\quad \forall t\in \mathbb{R}.
\]
\end{proof}
 \section{Conclusion} We have proved that any solution $x$ of ODE defined by a trigonometric polynomial field  can be approximated by a sequence functions  $t\mapsto \rho_k t+ \phi_k(t)$ where $(\rho_k)_k\subset \mathbb{Q}^n$ and converge to the rotation vector of $x$. The functions  $ \psi_k:\mathbb{R}\to\mathbb{R}^n$ are periodic on $t$  and uniformly bounded.
\makeatother

\appendix
\section*{Appendix. A}  
\begin{proof}[Proof of Lemma \ref{lem:cinfty}]
By hypothesis, for $q=0$ we have
\[
   \sum_{p \in 2\pi \mathbb{Z}^n/\mathbb{0}}\norm{c_{p }  } <+\infty,
\]
 the   series is normally convergent and we have
\[
\forall p\in 2\pi \mathbb{Z}^n:\ g(\omega z)\exp(- i  \langle z,p\rangle ) = \sum_{q \in 2\pi \mathbb{Z}^n}c_q\exp( i  \langle z,q-p\rangle ), \quad c_q \in \mathbb{C}^n,
\]
implies 
 \begin{align}
\label{equacprevided} a_p[g_\omega]&=\int_0^1\ldots \int_0^1g(\omega z)\exp(- i  \langle z,p\rangle )dz_1\ldots dz_n\\
\notag &=\sum_{q \in 2\pi \mathbb{Z}^n}c_q\int_0^1\ldots \int_0^1\exp( i  \langle z,q-p\rangle )dz_1\ldots dz_n=c_p.
\end{align}
Now, prove that  $g\in E_\omega(\mathbb{R}^n)$. Denote
\[
\theta_p(z):=\exp( i \frac{1}{\omega}\langle z,p\rangle ),\ \forall z\in \mathbb{R}^n,
\]
It is sufficient to prove that for every $q\geq1$ we have
\[
S_q:=\sum_{p\in  2\pi \mathbb{Z}^n/\mathbb{0} } \norm{c_{p}}  \norm{ d^q \theta_p}_\infty<+\infty,
\]
where $d^q  g$ is $q^{th}$ differential of $g$. The function $d^q  g$ is defined as
\[
d^q  g  =\sum_{p\in  2\pi \mathbb{Z}^n }  c_{p}   d^q \theta_p.
\]
 We have 
\begin{gather*}
  \forall p:=(p_j)_{j=1}^n\in 2\pi \mathbb{Z}^n:\\
  d^q \theta_p  =\sum_{k_1=1}^n\ldots \sum_{k_q=1}^n i^q \frac{1}{\omega^q}p_{k_1}\ldots p_{k_q}\exp( i \langle z,p\rangle ),
\end{gather*}
By consequence,
\[
 \norm{ d^q \theta_p}_\infty = \Big(\sum_{k=1}^n |\frac{p_k}{\omega}|\Big)^q\le n^q \Big(\frac{\norm{p}}{\omega}\Big)^q.
\]
Thanks to Equation \eqref{equacprevided}, we get
\begin{align*}
\forall q\geq1:\quad S_q&\le   n^q   \sum_{p\in  2\pi \mathbb{Z}^n/\mathbb{0} } \norm{c_{p}}  \Big(\frac{\norm{p}}{\omega}\Big)^q= n^q   \sum_{p\in  2\pi \mathbb{Z}^n/\mathbb{0} } \norm{a_p[g_\omega]}  \Big(\frac{\norm{p}}{\omega}\Big)^q \\
&=\frac{1}{2} n^q \norm{   g }_{\omega,q} <+\infty,
\end{align*}
which implies that  $g\in E_\omega(\mathbb{R}^n)$. 
\end{proof}
\section*{Appendix. B}  
\begin{proof}[ Proof of Lemma \ref{lem:equtow} ]
By Equation \eqref{lem:transfourier}, we have
\[
g(z)=\sum_{p \in  2\pi \mathbb{Z}^n }a_{p }[g_\omega] \exp(i \frac{1}{\omega}\langle z,p \rangle ).
\] 
Since $g(\mathbb{0})=\mathbb{0}, $ then
\[
a_{\mathbb{0}}[g_\omega]=-\sum_{p \in  2\pi \mathbb{Z}^n/\mathbb{0} }a_{p }[g_\omega],
\]
implies
\[
\norm{a_{\mathbb{0}}[g_\omega]}\le \sum_{p \in 2\pi \mathbb{Z}^n/\mathbb{0} }\norm{a_{p }[g_\omega]}=\frac{1}{2}\norm{   g}_{\omega,0}. 
\]
Since 
\[
\norm{   g}_\infty \le \sum_{p \in 2\pi \mathbb{Z}^n  }\norm{a_{p }[g_\omega] }.
\]
We  deduce that
\[
\norm{    g }_\infty\le  \norm{   g}_{\omega,0}.
\]
\end{proof}
\section*{Appendix. C}  
\begin{proof}[ Proof of Lemma \ref{lem:beta}]
Since $g: \mathbb{R}^n \to \mathbb{R}^n$  is a trigonometric polynomial, then it is $C^{\infty}$  and there exists $m\in \mathbb{N}$ such that 
\[
 g(z)=\sum_{p\in  2\pi \mathbb{Z}^n,\ \norm{p}\le m }a_{p}[g ] \exp(i \langle z,p\rangle ),\quad \forall z\in \mathbb{R}^n.
\] 
Implies,
\[ 
g_{\frac{1}{\omega}}(z)=g(\frac{z}{\omega})=\sum_{p\in  2\pi \mathbb{Z}^n,\ \norm{p}\le m }a_{p}[g ] \exp(i \frac{1}{\omega}\langle z,p\rangle ),
\]
Denote
\[
\theta_p(z):=\exp( i \frac{1}{\omega}\langle z,p\rangle ),\ \forall z\in \mathbb{R}^n,
\]
Denote $d^k  g$ is $k^{th}$ differential of $g$, which implies  
\begin{gather*}  
d^k  {\color{red} g_{\frac{1}{\omega}}(z)} =\sum_{p\in  2\pi \mathbb{Z}^n, \norm{p}\le m }  a_{p}[g ]   d^k \theta_p.
\end{gather*} 
 We have 
\begin{gather*}
\forall  p:=(p_j)_{j=1}^n\in  2\pi \mathbb{Z}^n:\\
 d^l \theta_p  =\sum_{s_1=1}^n\ldots \sum_{s_k=1}^n i^k \frac{p_{s_1}}{\omega}\ldots \frac{p_{s_k}}{\omega}\exp( i \frac{1}{\omega}\langle z,p\rangle  ,
\end{gather*}
then
\[
\forall p\in  2\pi \mathbb{Z}^n,\ \norm{p}\le m:\quad  \norm{ d^k \theta_p}_\infty = \Big(\sum_{s=1}^n |\frac{p_{s}}{\omega}|\Big)^k\le n^k \norm{\frac{m}{\omega}}^k.
\]
Since $\Big(g_{\frac{1}{\omega}}\Big)_\omega= g$., we obtain
\begin{align*}
\norm{   d^k   g_{\frac{1}{\omega}} }_{\omega,q}& \le 2  n^k  (\frac{m}{\omega})^{ k} \sum_{p\in 2\pi \mathbb{Z}^n, \ \norm{p}\le m }\norm{a_{p}[ g ] }   \norm{\frac{p}{\omega}}^{q}\\
&\le2  n^k  (\frac{m}{\omega})^{ k+q}\Big( \sum_{p\in 2\pi \mathbb{Z}^n, \ \norm{p}\le m }\norm{a_{p}[  g ] }  \Big).
\end{align*}
 It is sufficient to choose $\beta:=2 \max\{\sum_{p\in 2\pi \mathbb{Z}^n, \ \norm{p}\le m }\norm{a_{p}[   g ] }  , m\}$. 
\end{proof}
\section*{Appendix. D}  
\begin{proof}[Proof of Lemma \ref{lem:ineq} ]
Since $(h_j )_{j=1 }^k \subset  E_\omega(\mathbb{R}^n)$,   by Equation \eqref{lem:transfourier} we can write
\[
 h_j(z)=\sum_{p_{j}\in  2\pi \mathbb{Z}^n}a_{p_j}[h_{j,\omega}] \exp(i  \frac{1}{\omega}\langle z,p_j\rangle ).
\]
  By definition of the seminorm, 
\begin{equation}\label{equ:proofinequproduct}
\norm{    h_j}_{\omega,q}=2\sum_{p_{j}\in 2\pi \mathbb{Z}^n /\{\mathbb{0}\} }\norm{a_{p_j}[h_{j,\omega}]} \norm{\frac{p_j }{\omega} }^{q}.
\end{equation}
We have
\begin{align*}
 \Pi_{j=1}^k h_j(z)&=\sum_{p_{1},\ldots, p_{k}\in  2\pi \mathbb{Z}^n}  \Pi_{j=1}^k a_{p_j}[h_{j,\omega}] \exp( i\frac{1}{\omega} \langle z,\sum_{j=1}^k   p_j  \rangle )\\
 &=\sum_{v \in  2\pi \mathbb{Z}^n}  \sum_{\sum_{j=1}^k   p_j =v}\Pi_{j=1}^k a_{p_j}[h_{j,\omega}] \exp( i\frac{1}{\omega} \langle z,v \rangle ).
\end{align*}
Then
\begin{align*}
 \norm{    \Pi_{j=1}^k h_j }_{\omega,q}&=2 \sum_{v \in  2\pi \mathbb{Z}^n/\mathbb{0}} \norm{ \sum_{\sum_{j=1}^k   p_j =v}\Pi_{j=1}^k a_{p_j}[h_{j,\omega}]} \norm{\frac{v}{\omega}}^q\\
 &\le 2 \sum_{v \in  2\pi \mathbb{Z}^n } \norm{ \sum_{\sum_{j=1}^k   p_j =v}\Pi_{j=1}^k a_{p_j}[h_{j,\omega}]} \norm{\frac{v}{\omega}}^q.
\end{align*}
Using the triangular inequality, we obtain
\begin{align*}
 \norm{    \Pi_{j=1}^k h_j }_{\omega,q}&\le 2 \sum_{p_{1},\ldots, p_{k}\in 2\pi \mathbb{Z}^n }  \norm{\Pi_{j=1}^k  a_{p_j}[h_{j,\omega}] }\norm{   \frac{1}{\omega} \sum_{j=1}^k {p_j} }^{q}.
\end{align*}
Since 
\[
 \norm{    \sum_{j=1}^k {p_j} }  \le    k  \Pi_{\norm{ p_j }\neq 0}  \norm{ p_j }, \quad \forall p_j,\in 2\pi \mathbb{Z}^n,
\]
then
\[
 \norm{    \frac{1}{\omega}  \sum_{j=1}^k {p_j} }^q  \le     k^q \omega^{(k-1)q} \Pi_{\norm{ p_j }\neq 0}  \norm{ \frac{p_j }{\omega} }^q, \quad \forall p_j,\in 2\pi \mathbb{Z}^n.
\]
We deduce that
\begin{align*}
&\norm{   \Pi_{j=1}^k h_j}_{\omega,q}\le  2 k^{q}\omega^{(k-1)q} \sum_{p_{1},\ldots, p_{k}\in 2\pi \mathbb{Z}^n }   \Big( \Pi_{\norm{  p_j }\neq 0}  \norm{  \frac{p_j }{\omega} }^{q}\Big)\Big( \Pi_{j=1}^k \norm{  a_{p_j}[h_{j,\omega}] }\Big) \\\\
&=2 k^{q}\omega^{(k-1)q} \sum_{p_{1},\ldots, p_{k}\in  2\pi \mathbb{Z}^n}  \Big( \Pi_{ \norm{ p_j }= 0 }\norm{  a_{p_j}[h_{j,\omega}] } \Big)  \Pi_{\norm{  p_j }\neq 0}  \Big[\norm{  \frac{p_j }{\omega} }^{q}   \norm{  a_{p_j}[h_{j,\omega}] }\Big].
\end{align*} 
By Equation \eqref{equ:proofinequproduct},
\begin{align*}
\norm{\Pi_{j=1}^k h_j }_{\omega,q}&\le  \frac{1}{2^{ k-1 }}k^{q} \omega ^{(k-1)q}    \Pi_{j=1}^k\Big[ 2\Big( \norm{a_0[h_{j,\omega}]}+\sum_{p_{j}\in 2\pi \mathbb{Z}^n/\{\mathbb{0}\}}  \norm{ \frac{p_j }{\omega} }^{q} \norm{a_{p_j}[h_{j,\omega}]}\Big)\Big].
\end{align*} 
\end{proof}
\end{document}